\documentclass[11pt]{article}
\usepackage{hyperref}
\usepackage{graphicx,amssymb,upref}
\usepackage{latexsym}
\usepackage{amsthm}
\usepackage{mathtools}
\usepackage{amsfonts, epsfig, amsmath, color, amscd}
\usepackage{pb-diagram, pb-xy}
\newtheorem{conjecture}{Conjecture}

\theoremstyle{plain}

\newtheorem{theorem}{Theorem}
\newtheorem{lemma}{Lemma}

\title{On a lower bound for the energy functional on a family of Hamiltonian minimal Lagrangian tori in $\mathbb{C}P^2$} 

\author{A. A. Kazhymurat}  
\date{29 September 2017}

\begin{document} 
 \maketitle
\begin{abstract}  
      We study the energy functional on the set of Lagrangian tori in $\mathbb{C}P^{2}$. We prove that the value of the energy functional on a certain family of Hamiltonian minimal Lagrangian tori in  $\mathbb{C}P^{2}$ is strictly larger than energy of the Clifford torus.
\end{abstract} 

%%%% **** The text of the paper starts here **** %%%%

\section{Introduction}
As remarked in \cite{operator}, one can naturally associate a 2D periodic Schr\"odinger operator with every Lagrangian torus in $\mathbb{C}P^{2}$. More precisely, any Lagrangian torus $\Sigma \subset \mathbb{C}P^{2}$ with induced metric
\begin{equation}
\label{metric}
ds^2=2 e^{v(x, y)}(dx^2+dy^2)
\end{equation} is the image of the composition of mappings 
$$
r: \mathbb{R}^{2} \rightarrow S^{5} \xrightarrow{\mathcal{H}} \mathbb{C}P^2,
$$
where $r$ is a horizontal lift, $\mathcal{H}$ is the Hopf projection. The vector-function $r$ satisfies the Schr\"odinger equation
$$
Lr=0,\qquad L=(\partial_{x}-\frac{i \beta_x}{2})^{2}+(\partial_{y}-\frac{i \beta_y}{2})^{2}+V(x, y),\qquad V=4 e^{v}+\frac{1}{4}(\beta_x^2+\beta_y^2)+\frac{i}{2}\Delta \beta,
$$
where $\beta$ is the Lagrangian angle (see the definition below).

The existence of operator $L$ allows us to introduce the energy functional $E$ on the set of Lagrangian tori in $\mathbb{C}P^{2}$ (see \cite{energy})
$$
E(\Sigma)=\frac{1}{2}\int_{\Sigma}V \, dx \wedge dy.
$$
As shown in \cite{energy} the energy functional admits following geometric interpretation
$$
E(\Sigma)=A(\Sigma)+\frac{1}{8}W(\Sigma), \qquad A(\Sigma)=\int_{\Sigma}d \sigma,\qquad W(\Sigma)=\int_{\Sigma} |H|^{2} \, d\sigma,
$$
where $d\sigma=2 e^{v} dx \wedge dy$ is the induced area element, $H$ is the mean curvature vector field.

For the Clifford torus $\Sigma_{Cl}$ which is defined by the vector-function
$$
r(x, y)=\bigl(\frac{1}{\sqrt{3}} e^{2 \pi i x},\frac{1}{\sqrt{3}} e^{2 \pi i(-\frac{1}{2}x+\frac{\sqrt{3}y}{2})},\frac{1}{\sqrt{3}}  e^{2 \pi i (-\frac{1}{2}x-\frac{\sqrt{3}y}{2})}\bigr),
$$
energy equals
$$
E(\Sigma_{Cl})=\frac{4 \pi^{2}}{3 \sqrt{3}}.
$$
Following conjecture was proposed in \cite{energy}.
\begin{conjecture}
{ \it The minimum of the energy functional is attained on the Clifford torus.}
\end{conjecture}
Conjecture 1 has been verified for two families of Hamiltonian minimal Lagrangian tori: for homogeneous tori and for tori constructed in \cite{ToriExamples}.

A homogeneous torus $\Sigma_{r_1, r_2, r_3} \subset \mathbb{C}P^{2}, r_1^2+r_2^2+r_3^2=1,  r_i>0$ is defined by the vector-function
$$
r(x, y)= \bigl(r_1 e^{2 \pi i x}, r_2 e^{2 \pi i (a_1 x+b_1 y)},r_3 e^{2 \pi i (a_2 x+b_2 y)}\bigr),
$$
with some restrictions on $a_i, b_i$. Following inequality holds
$$
E(\Sigma_{r_1, r_2, r_3})=\frac{\pi^{2}(1-r_1^2)(1-r_2^2)(1-r_3^2)}{2 r_1 r_2 r_3} \geqslant \frac{4\pi^{2}}{3 \sqrt{3}},
$$
and equality is attained only for the Clifford torus.

The second family of tori $\Sigma_{m, n, k} \subset \mathbb{C}P^{2}, m, n, k \in \mathbb{Z}, m \geqslant n>0, k<0$ is of form $\mathcal{H}(\tilde{\Sigma}_{m, n, k})$ where
$$
\tilde{\Sigma}_{m, n, k}=\left\{(u_1 e^{2 \pi i m y}, u_2 e^{2 \pi i n y}, u_3 e^{2 \pi i k y})\right\} \subset S^{5},
$$
the numbers $u_1, u_2, u_3$ satisfy the equation
$$
u_1^2+u_2^2+u_3^2=1, \qquad m u_1^2+n u_2^2+k u_3^2=0.
$$
The parameters $m, n, k$ should be chosen so that the involution
$$
(u_1, u_2, u_3) \longrightarrow (u_1 \mathrm{cos}(m \pi), u_2 \mathrm{cos}(n \pi), u_3 \mathrm{cos} (k \pi))
$$
on the surface $m u_1^2+n u_2^2+k u_3^2=0$ preserves its orientation (otherwise $\mathcal{H}(\tilde{\Sigma}_{m, n, k})$ is homeomorphic to Klein bottle, see \cite{ToriExamples}). Following inequality is proved in \cite{energy} 
$$
E(\Sigma_{m, n, k})>E(\Sigma_{Cl}).
$$
In the case of minimal Lagrangian tori the function $v(x, y)$ satisfies the Tzizeica equation (see \cite{MironovTorus}). Smooth periodic solutions of this equation are finite-gap, i.e. can be expressed in terms of the theta-function on the Jacobian variety of the spectral curve. The results of \cite{hask} imply the conjecture for minimal Lagrangian tori corresponding to spectral curve of sufficiently high genus.

One should note that for embedded Lagrangian tori with non-trivial Floer cohomology one can derive lower bounds for the area functional. For instance, following inequality holds for any Lagrangian torus $\Sigma$ Hamiltonian isotopic to the Clifford torus \cite{gold}
$$
A(\Sigma) \geqslant \frac{3}{\pi} E(\Sigma_{Cl}).
$$
It is unclear at present whether one can derive symplectic-topological bounds for the Willmore functional. This seems to be related to the question whether every monotone Lagrangian torus in $\mathbb{C}P^2$ is Hamiltonian isotopic to a minimal torus.

The aim of the present work is to verify the conjecture 1 for the family of Hamiltonian minimal Lagrangian tori constructed in \cite{MironovTorus} (also see \cite{MaTorus}).

Let $\alpha_1, \alpha_2, \alpha_3 \in \mathbb{Z}, b=-\alpha_1-\alpha_2-\alpha_3, c=\alpha_1\alpha_2+\alpha_1\alpha_3+\alpha_2\alpha_3, c_1=-\alpha_1\alpha_2\alpha_3$, $a_1>a_2>0$ be some real numbers satisfying the inequalities \eqref{ssylka15}, \eqref{ssylka16} (see below). Following theorem has been proved in \cite{MironovTorus}.
\begin{theorem}
The mapping $\psi: \mathbb{R}^{2} \rightarrow \mathbb{C}P^{2}$ defined by the formula
$$
\psi(x, y)=\bigl(F_1(x)e^{i(G_1(x)+\alpha_1 y)}:F_2(x)e^{i(G_2(x)+\alpha_2 y)}:F_3(x)e^{i(G_3(x)+\alpha_3 y)}\bigr),
$$
is a conformal Hamiltonian minimal Lagrangian immersion, where
$$
F_i=\sqrt{\frac{2 e^{v}+\alpha_{i+1}\alpha_{i+2}}{(\alpha_i-\alpha_{i+1})(\alpha_{i}-\alpha_{i+2})}}, \qquad G_i=\alpha_i\int_{0}^{x}\frac{c_2-a e^{v}}{2 \alpha_i e^{v}-c_1}\, dz,
$$
\begin{equation}
2 e^{v(x)}=a_1\left(1-\frac{a_1-a_2}{a_1}\mathrm{sn}^{2}\left(x \sqrt{a_1+a_3}, \frac{a_1-a_2}{a_1+a_3}\right)\right)
\end{equation}
(index $i$ runs modulo 3), $\mathrm{sn}(x)$ is the Jacobi's elliptic function, $c_2$ is a real root of \eqref{ssylkac_2}, $a_3=\frac{c_1^2+c_2^2}{a_1 a_2}$.
\end{theorem}

Moreover, if the rationality constraints \eqref{rational} are met, $\psi$ is a doubly periodic mapping and the image of the plane is a Hamiltonian minimal Lagrangian torus $\Sigma_{M} \subset \mathbb{C}P^{2}$.

The principal result of the present work is following theorem.
\begin{theorem}
\label{theorema}
The inequality
$$
E(\Sigma_M)>E(\Sigma_{Cl})
$$
holds if $\alpha_1-\alpha_3, \alpha_2-\alpha_3$ are relatively prime.
\end{theorem}
The theorem \ref{theorema} thus confirms the conjecture 1.

\begin{section}{The proof of the theorem \ref{theorema}}
Lagrangianity of $\Sigma$, horizontality of the mapping $r:\mathbb{R}^{2} \rightarrow S^{5}$ and the form of the induced metric \eqref{metric} imply
$$
R=\begin{pmatrix}
 r\\
 \frac{r_x}{|r_x|}\\
 \frac{r_y}{|r_y|}\\
 \end{pmatrix}\in U(3).
$$
The Lagrangian angle $\beta(x, y)$ is defined by the equation $e^{i \beta}=\mathrm{det} R$.
The mean curvature vector field can be expressed in terms of the Lagrangian angle $H=J \nabla \beta$ where $J$ is the complex structure on $\mathbb{C}P^{2}$. For minimal tori $\beta=\mathrm{const}$. As demonstrated in \cite{oh} in the case of Hamiltonian minimal tori $\beta$ is a linear function in the conformal coordinates $x, y$.

Let us consider the Hamiltonian minimal immersion $\psi$ \cite{MironovTorus} defined in the theorem \ref{theorema}.

The equation
$$
(a_1-a_2)^{2} x^4+2(a_1^3 a_2^2+a_1^2 a_2^3+(a_1^2 a_2+a_1 a_2^2)b c_1+(a_1^2+a_2^2)c_1^2+2 a_1^2 a_2^2 c) x^2+
$$
\begin{equation}
\label{ssylkac_2}
+((a_1+a_2)c_1^2-a_1^2 a_2^2+a_1 a_2 b c_1)^2=0.
\end{equation}
has a real root $x=c_2$ iff following inequalities are satisfied
\begin{equation}
P=a_1^3 a_2^2+a_1^2 a_2^3+(a_1^2 a_2+a_1 a_2^2)b c_1+(a_1^2+a_2^2)c_1^2+2 a_1^2 a_2^2 c \leqslant 0, 
\label{ssylka15}
\end{equation}
\begin{equation}
P^2-(a_1-a_2)^2 ((a_1+a_2)c_1^2-a_1^2 a_2^2+a_1 a_2 b c_1)^2 \geqslant 0.
\label{ssylka16}
\end{equation}
Recall that $\mathrm{sn}(u, k)=\mathrm{sin} \, \theta$ where
\begin{equation}
\label{elliptic}
u(\theta)=\int_{0}^{\theta} \frac{d\phi}{\sqrt{1-k^{2}\mathrm{sin}^{2}\phi}}.
\end{equation}
The function $\mathrm{sn}^2(u)$ is periodic with period $2u(\frac{\pi}{2})$ (see, for instance, \cite{Akhiezer}). Therefore $v(x)$ has period
\begin{equation}
T=\frac{2 u\left(\frac{\pi}{2}\right)}{\sqrt{a_1+a_3}}.
\label{periodT}
\end{equation}
Further we assume that $(\alpha_1-\alpha_3, \alpha_2-\alpha_3)=1$. 

The immersion $\psi:\mathbb{R}^{2}\rightarrow \mathbb{C}P^{2}$ is doubly periodic if there exists $\tau \in \mathbb{R}$ such that
\begin{equation}
\label{rational}
\lambda_1=\frac{G_1(T)-G_3(T)+(\alpha_1-\alpha_3)\tau}{2 \pi}, \qquad \lambda_2=\frac{G_2(T)-G_3(T)+(\alpha_2-\alpha_3)\tau}{2 \pi} \in \mathbb{Q}.
\end{equation}
Then the vectors of period can be expressed as follows
$$
e_1=(0, 2\pi), \qquad e_2=N(T, \tau),
$$
where $N$ is some natural number.
If the condition \eqref{rational} is met, $\Sigma_M \subset \mathbb{C}P^{2}$ is an immersed torus with Lagrangian angle $\beta=ax+by$ where
$$
a=\frac{b c_1+a_1 a_3+a_2 a_3-a_1 a_2}{c_2}.
$$
Following equality holds
$$
|H|^{2}=\frac{1}{2}e^{-v}(a^2+b^2).
$$
Let us find lower bounds for $W(\Sigma_M)$ and $A(\Sigma_M)$.

Using \eqref{periodT} and $a_3>0$ we arrive at the inequalities
$$
u(\frac{\pi}{2})>\frac{\pi}{2}, \qquad T>\frac{\pi}{\sqrt{a_1+a_3}}.
$$
Thus
$$
W(\Sigma_{M})=\int_{\Sigma_{M}} \lvert H \rvert^{2}\, d\sigma=\int_{\Lambda}\frac{1}{2}e^{-v}(a^{2}+b^{2}) 2 e^{v}\, dx \wedge dy=2\pi N T (a^2+b^2).
$$
Therefore, following lower bound for $W(\Sigma_{M})$ holds
\begin{equation}
\label{willmoreestimation}
W(\Sigma_{M})>2\pi^{2}\frac{a^{2}+b^2}{\sqrt{a_1+a_3}}.
\end{equation}

Following lemma provides a lower bound for $A(\Sigma_M)$.
\begin{lemma}
\label{areaestimation}
The inequality
$$
A(\Sigma_{M})>\pi^{2}\frac{a_1+a_2}{\sqrt{a_1+a_3}}
$$
is true.
\end{lemma}
Proof of the lemma \ref{areaestimation}. We have
$$
A(\Sigma_{M})=\int_{\Sigma_{M}}\, d\sigma=\int_{\Lambda}2 e^{v(x)} \, dx\wedge dy=2\pi\int_{0}^{NT}2 e^{v(x)} dx \geqslant 2\pi \int_{0}^{T}2e^{v(x)}\, dx=
$$
$$
=2\pi\int_{0}^{T}a_1\left(1-\frac{a_1-a_2}{a_1}\mathrm{sn}^{2}\left(x\sqrt{a_1+a_3},\frac{a_1-a_2}{a_1+a_3}\right)\right)\, dx=
$$
$$
=\frac{2 \pi a_1}{\sqrt{a_1+a_3}}\int_{0}^{2u(\frac{\pi}{2})}\left(1-\frac{a_1-a_2}{a_1}\mathrm{sn}^2\left(u, \frac{a_1-a_2}{a_1+a_3}\right)\right)\, du.
$$

Using \eqref{elliptic} we arrive at
$$
\int_{0}^{T}2 e^{v(x)}\, dx=\frac{a_1}{\sqrt{a_1+a_3}}\int_{0}^{\pi}\frac{1-\frac{a_1-a_2}{a_1}\mathrm{sin}^{2}\theta}{\sqrt{1-\left(\frac{a_1-a_2}{a_1+a_3}\right)^{2}\mathrm{sin}^{2}\theta}}\, d\theta.
$$
As $0<\frac{a_1-a_2}{a_1+a_3}<1$, following estimate is true
$$
\int_{0}^{T}2 e^{v(x)} dx>\frac{a_1}{\sqrt{a_1+a_3}}\int_{0}^{\pi}\left(1-\frac{a_1-a_2}{a_1}\mathrm{sin}^{2}\theta\right)\, d\theta=\frac{\pi(a_1+a_2)}{2\sqrt{a_1+a_3}}.
$$
Lemma \ref{areaestimation} is proved.

The inequalities \eqref{ssylka15}, \eqref{ssylka16} are invariant under simultaneous change of sign  $\alpha_{1},\alpha_{2},\alpha_{3}$ and their permutations. If $\alpha_1, \alpha_2,\alpha_3$ are all of the same sign, the inequality \eqref{ssylka15} has no positive solutions. Therefore we assume without loss of generality that $\alpha_1\geqslant\alpha_2\geqslant 0\geqslant\alpha_3$.

\begin{lemma}  
\label{lemmaoninequality}
If $\alpha_1\geqslant\alpha_2\geqslant0\geqslant\alpha_3$ и $a_1>a_2>0$, the inequalities (\ref{ssylka15}) and (\ref{ssylka16}) are satisfied simultaneously iff
\begin{equation}
    \label{square}
-\alpha_2\alpha_3\leqslant
a_2<a_1
\leqslant
-\alpha_1\alpha_3.
\end{equation}
\end{lemma}

Proof of the lemma \ref{lemmaoninequality}. 
Denote
$$
Q(x)=-(x+\alpha_1 \alpha_2)(x+\alpha_1 \alpha_3)(x+\alpha_2\alpha_3).
$$
Then \eqref{ssylkac_2} assumes the form
$$
(a_1-a_2)^2\left(x^2-\left(\frac{a_1\sqrt{Q(a_2)}-a_2\sqrt{Q(a_1)}}{a_1-a_2}\right)^2\right)\left(x^2-\left(\frac{a_1 \sqrt{Q(a_2)}+a_2 \sqrt{Q(a_1)}}{a_1-a_2}\right)^2\right)=0.
$$
This equation has a positive root iff $Q(a_1)\geqslant 0, Q(a_2) \geqslant 0$. This is equivalent to $
-\alpha_2\alpha_3\leqslant
a_2<a_1
\leqslant
-\alpha_1\alpha_3.
$
Lemma \ref{lemmaoninequality} is proved.

It follows from the proof of the lemma \ref{lemmaoninequality} that if $\alpha_3=0$ or $\alpha_1=\alpha_2$ inequalities \eqref{ssylka15}, \eqref{ssylka16} are not satisfied for $a_1>a_2$. Therefore we assume without loss of generality
\begin{equation}
    \label{alpha}
\alpha_1>\alpha_2 \geqslant 0>\alpha_3.
\end{equation}
The inequality (\ref{willmoreestimation}) and lemma \ref{areaestimation} imply 
$$
E(\Sigma_{M})>\pi^{2}\frac{a_1+a_2+\frac{a^{2}+b^2}{4}}{\sqrt{a_1+a_3}}.
$$

Let us prove $E(\Sigma_M)>E(\Sigma_{Cl})$. We will consider two cases: $\alpha_2 >0$ and $\alpha_2=0$.

Assume $\alpha_2 >0$.

If $(a_1+a_2)a_3 \geqslant \frac{7}{4}(a_1 a_2-b c_1)$ then
   $$
   a^{2}=\frac{((a_1+a_2)a_3-(a_1 a_2-b c_1))^2}{c_2^2} \geqslant \frac{9}{49} (a_1+a_2)^2 \frac{a_3^2}{c_2^2}=\frac{9}{49} (a_1+a_2)^2 \frac{a_3}{a_1 a_2}\frac{c_1^2+c_2^2}{c_2^2}\geqslant \frac{9}{49}(a_1+a_2)^2\frac{a_3}{a_1 a_2}.
   $$ 
   As $a_1>a_2 \geqslant 1$ и $(a_1+a_2)^{2} > 4 a_1 a_2$ we have
   $$
   E(\Sigma_{M})>\pi^{2}\frac{a_1+a_2+\frac{9(a_1+a_2)^2 a_3}{196 a_1 a_2}}{\sqrt{a_1+a_3}} >\pi^{2} \frac{a_1+\frac{9a_3}{49}}{\sqrt{a_1+a_3}}=\pi^{2} \sqrt{a_1} \frac{1+\frac{9a_3}{49a_1}}{\sqrt{1+\frac{a_3}{a_1}}}>\pi^{2}\frac{1+\frac{9a_3}{49a_1}}{\sqrt{1+\frac{a_3}{a_1}}}.
   $$
Note that for positive $x$ we have $
\frac{1+\frac{9x}{49}}{\sqrt{1+x}}>\frac{4}{3\sqrt{3}}$ holds.
Consequently, $E(\Sigma_M)>E(\Sigma_{Cl})$.

Now consider the case
$$
(a_1+a_2)a_3 < \frac{7}{4}(a_1 a_2-b c_1).
$$ We analyse two cases: $\alpha_1>-\frac{3}{2}\alpha_2\alpha_3$ and $\alpha_1 \leqslant -\frac{3}{2}\alpha_2\alpha_3$.
   
     If $\alpha_1> -\frac{3}{2} \alpha_2\alpha_3$ then
   $$
  \alpha_1 < -3 b=3(\alpha_1+\alpha_2+\alpha_3),
   $$
   as $\alpha_1>-\frac{3}{2}(\alpha_2+\alpha_3)$. From \eqref{square}
   $$
   -\frac{b c_1}{a_1+a_2}=\frac{b \alpha_1 \alpha_2 \alpha_3}{a_1+a_2} < \frac{b (3 b) \alpha_2 \alpha_3 }{2 \alpha_2 \alpha_3}=\frac{3}{2} b^2.
   $$
   Hence
   $$
   E(\Sigma_{M})>\pi^{2}\frac{a_1+a_2+\frac{b^2}{4}}{\sqrt{a_1+a_3}}>\pi^{2}\frac{a_1+a_2+\frac{b^2}{4}}{\sqrt{a_1+\frac{7}{4}\frac{a_1a_2}{a_1+a_2}-\frac{7}{4}\frac{b c_1}{a_1+a_2}}}>\pi^2 \frac{a_1+a_2+\frac{b^2}{4}}{\sqrt{a_1+\frac{7}{4} a_2+\frac{21}{8} b^2}}> 
   $$
   $$
   >\pi^2 \frac{a_1+a_2+\frac{b^2}{4}}{\sqrt{\frac{7}{4}a_1+\frac{7}{4} a_2+\frac{21}{8} b^2}}
   = \pi^2 \sqrt{\frac{4(a_1+a_2)}{7}}\frac{1+\frac{b^2}{4(a_1+a_2)}}{\sqrt{1+\frac{3}{2}\frac{b^2}{a_1+a_2}}}> \pi^2 \sqrt{\frac{8}{7}} \frac{1+\frac{b^2}{4(a_1+a_2)}}{\sqrt{1+\frac{3}{2}\frac{b^2}{a_1+a_2}}}>E(\Sigma_{Cl}).
   $$
   The last inequality can be seen by considering the function $f(x)=\sqrt{\frac{8}{7}}\frac{1+\frac{x}{4}}{\sqrt{1+\frac{3}{2}x}}$ for $x>0$.

    If $\alpha_1\leqslant -\frac{3}{2} \alpha_2\alpha_3$, the inequalities \eqref{square} and \eqref{alpha} imply
   $$
   -b c_1\leqslant - 2\alpha_1^2  \alpha_2 \alpha_3 < \frac{9}{2} a_1 a_2^2.
   $$
   Therefore
   $$
   E(\Sigma_M)>\pi^{2}\frac{a_1+a_2}{\sqrt{a_1+\frac{7}{4} \frac{a_1 a_2-b c_1}{a_1+a_2}}}=\pi^2 \frac{(a_1+a_2)\sqrt{a_1+a_2}}{\sqrt{a_1(a_1+a_2)+\frac{7}{4}a_1 a_2-\frac{7}{4}b c_1}}>
   $$
   $$
   >\pi^{2}\frac{(a_1+a_2)\sqrt{a_1+a_2}}{\sqrt{a_1^2+\frac{11}{4}a_1 a_2+\frac{63}{8} a_1 a_2^2}}
   >\pi^{2}\frac{(a_1+a_2)\sqrt{a_1+a_2}}{\sqrt{a_1^3+\frac{11}{4}a_1^2 a_2+\frac{63}{8}  a_1 a_2^2}}=\pi^{2} \frac{(1+\frac{a_2}{a_1})\sqrt{1+\frac{a_2}{a_1}}}{\sqrt{1+\frac{11}{4}\frac{a_2}{a_1}+\frac{63}{8}\frac{a_2^2}{a_1^2}}}>E(\Sigma_{Cl}).
   $$

Let us consider the case $\alpha_2=0$. Introduce $p=-\alpha_1 \alpha _3$, $x=\frac{a_{1}}{p},y=\frac{a_{2}}{p}$. Note that $0<y<x\leqslant 1$ due to \eqref{alpha}. Then inequalities \eqref{ssylka15}, \eqref{ssylka16} assume following form
$$
p^5 x^2 y^2 (x+y-2) \leqslant 0, \qquad 4 p^{10} x^4 y^4 (1-x) (1-y) \geqslant 0.
$$
The equation \eqref{ssylkac_2} implies
\begin{equation}
\label{degenerate}
c_{2}^2=p^{3} x^2 y^2 \frac{2-x-y \pm \sqrt{(2-x-y)^2-(x-y)^2}}{(x-y)^2}.
\end{equation}
As $2-x-y>0$ we have $\sqrt{(2-x-y)^2-(x-y)^2}=(2-x-y)\sqrt{1-\frac{(x-y)^2}{(2-x-y)^2}}.$
Note that by Bernoulli inequality $$
1-\frac{(x-y)^2}{(2-x-y)^2}\leqslant \sqrt{1-\frac{(x-y)^2}{(2-x-y)^2}} \leqslant 1-\frac{(x-y)^2}{2(2-x-y)^2}.
$$ 
Consequently,
\begin{equation}
\label{squareroot}
2-x-y-\frac{(x-y)^2}{2-x-y} \leqslant \sqrt{(2-x-y)^2-(x-y)^2} \leqslant 2-x-y-\frac{(x-y)^2}{2(2-x-y)}.
\end{equation}
Consider two cases: sign '+' and '-' in \eqref{degenerate}. For the '-' sign \eqref{degenerate} and \eqref{squareroot} imply the inequalities
$$p^{3} \frac{x^{2}y^{2}}{2(2-x-y)} \leqslant c_{2}^{2}\leqslant p^{3}\frac{x^{2}y^{2}}{2-x-y}.$$ 
As $c_1=0$ we have following bound for $a_3$
$$
a_3=\frac{c_2^2}{a_1 a_2}, \qquad
p\frac{xy}{2(2-x-y)}\leqslant a_{3} \leqslant p\frac{xy}{2-x-y}.
$$
These estimates and lemma \ref{areaestimation} imply
$$
A(\Sigma_{M}) \geqslant \pi^{2} \sqrt{p} \frac{x+y}{\sqrt{x+\frac{xy}{2-x-y}}}.
$$
Following inequality holds
$$
a=\frac{(a_1+a_2)a_3-a_1 a_2}{c_2} \geqslant \frac{(xp+yp)p\frac{xy}{2(2-x-y)}-xyp^2}{c_2} \geqslant \sqrt{p}\left(\frac{x+y}{2(2-x-y)}-1\right)\sqrt{2-x-y}.
$$
The estimate \eqref{willmoreestimation} implies
$$
W(\Sigma_{M}) \geqslant 2\pi^2 \frac{a^2}{\sqrt{a_1+a_3}}\geqslant 2 \pi^{2} \sqrt{p} \left(\frac{x+y}{2(2-x-y)}-1\right)^{2}\frac{2-x-y}{\sqrt{x+\frac{xy}{2-x-y}}}.
$$
Henceforth
$$
E(\Sigma_{M}) \geqslant \pi^{2} \sqrt{p}\left(\frac{x+y}{\sqrt{x+\frac{xy}{2-x-y}}}+\frac{1}{4}\left(\frac{x+y}{2(2-x-y)}-1\right)^{2}\frac{2-x-y}{\sqrt{x+\frac{xy}{2-x-y}}}\right).
$$
As $p\geqslant 1$ we have
$$
E(\Sigma_M) \geqslant \pi^2 B_1(x, y), \qquad B_1(x, y)=\frac{16-7x^2+8x-14yx+8y-7y^2}{16\sqrt{(2-x)(2-x-y)x}}.
$$
\begin{lemma}
\label{B_1}
If $0<y<x\leqslant 1$, then $B_1(x, y)>1$.
\end{lemma}
Proof of the lemma \ref{B_1}. One can check by direct computation that there are no critical points $\partial_x B_1=\partial_y B_1=0$ inside the triangle $0<y<x\leqslant 1$ while on the boundary of the triangle $B_1(x, y)>1$ holds. Lemma \ref{B_1} is proved.

Therefore, $E(\Sigma_M)>E(\Sigma_{Cl})$ holds for the '-' sign in \eqref{degenerate}.

For the '+' sign in \eqref{degenerate}  \eqref{squareroot} implies the inequalities
$$
p^{3} f(x, y) \leqslant c_{2}^{2} \leqslant p^{3} g(x,y), 
$$
where
$$
f(x,y)=x^{2}y^{2}\frac{2(2-x-y)-\frac{(x-y)^{2}}{2-x-y}}{(x-y)^{2}}, \qquad g(x,y)=x^{2}y^{2}\frac{2(2-x-y)-\frac{(x-y)^{2}}{2(2-x-y)}}{(x-y)^{2}}.
$$
Analogously one establishes the inequalities
$$
p \frac{f(x,y)}{xy}\leqslant a_{3} \leqslant p \frac{g(x,y)}{xy},
$$
$$
a \geqslant \sqrt{p}\frac{(x+y)\frac{f(x,y)}{xy}-xy}{\sqrt{g(x,y)}}.
$$
The inequality \eqref{willmoreestimation} and lemma \ref{areaestimation} imply
$$
A(\Sigma_M) \geqslant \pi^{2} \sqrt{p} \frac{x+y}{\sqrt{x+\frac{g(x,y)}{xy}}},
$$
$$
W(\Sigma_M)\geqslant 2\pi^2 \frac{a^2}{\sqrt{a_1+a_3}} \geqslant 2 \pi^{2} \sqrt{p} \frac{((x+y)\frac{f(x,y)}{xy}-xy)^2}{g(x, y)\sqrt{x+\frac{g(x,y)}{xy}}},
$$
$$
E(\Sigma_M) \geqslant
    \pi^{2} \sqrt{p}\frac{x+y+\frac{1}{4}\frac{((x+y)\frac{f(x,y)}{xy}-xy)^2}{g(x,y)}}{\sqrt{x+\frac{g(x,y)}{xy}}} \geqslant \pi^{2}  B_2(x,y),$$ 
    где $$
    B_2(x, y)= \frac{x+y+\frac{1}{4}\frac{((x+y)\frac{f(x,y)}{xy}-xy)^2}{g(x,y)}}{\sqrt{x+\frac{g(x,y)}{xy}}}.
    $$
The following lemma is established similarly to the lemma \ref{B_1}.
\begin{lemma}
If $0<y<x\leqslant 1$, then $B_1(x, y)>0.9$.
\end{lemma} 
This finishes the proof of the theorem \ref{theorema}.

\end{section}
\bibliographystyle{amsplain}

\begin{thebibliography}{10}
\bibitem {operator} A.E.~Mironov. The Novikov-Veselov hierarchy of equations and integrable deformations of
minimal Lagrangian tori in $\mathbb{C}P^{2}$. (Russian) Sib. Elektron. Mat. Izv., 2004, $\textbf{1}$, 38--46.

\bibitem {energy} H.~Ma, A.E.~Mironov, D.~Zuo. Energy functional for Lagrangian tori in $\mathbb{C}P^{2}$, \href{https://arxiv.org/abs/1701.07211}{arXiv:1701.07211}.


\bibitem {ToriExamples}  A.E.~Mironov. 	New examples of Hamilton-minimal and minimal Lagrangian submanifolds in $\mathbb{C}^{n}$ and $\mathbb{C}P^{n}$. Sbornik: Mathematics, 2004, $\textbf{195}$:1, 85--96.

\bibitem {gold} E.~Goldstein. Some estimates related to Oh's conjecture for the Clifford tori in $\mathbb{C}P^{n}$, \href{https://arxiv.org/abs/math/0311460}{arXiv:math/0311460}.

\bibitem{MironovTorus} A.E.~Mironov. On Hamiltonian-minimal Lagrangian tori in $\mathbb{C}P^{2}$. Sib. Math. J., 2003 $\textbf{44}$:6,  1324--1328.

\bibitem {hask} M.~Haskins. The geometric complexity of special Lagrangian $T^{2}-$cones. Invent. Math., 2004, \textbf{157}:1, 11--70.

\bibitem {MaTorus} 
H.~Ma, M.~Schmies. Examples of Hamiltonian stationary Lagrangian tori in $\mathbb{C}P^{2}$. Geom.
Dedicata, 2006, $\textbf{118}$, 173--183.

\bibitem {oh} Y.~Oh. Volume minimization of Lagrangian submanifolds under Hamiltonian deformations,
Math. Z. 1993, \textbf{212}, 175--192.

\bibitem {Akhiezer} N.I.~Akhiezer. Elements of the theory of elliptic functions. Transl. Math. Monogr., 1990, \textbf{79},  p. 208.
\end{thebibliography}

\bigskip

NIS of Physics and Mathematics

Zhamakaev St, Almaty 55000, Kazakhstan

\href{mailto:akkazhymurat@gmail.com}{akkazhymurat@gmail.com}
\end{document}